\documentclass[12pt]{amsart}
\usepackage{amssymb,latexsym,amscd,amsfonts,pb-diagram}
\usepackage{mathrsfs}
\usepackage[all]{xy}
\usepackage{color}
\newcommand{\incto}{\ar@{^{(}->}}
\newcommand{\into}{\ar@{>->}}
\newcommand{\onto}{\ar@{->>}}
\let\epsilon=\varepsilon

\newtheorem{Theorem}{Theorem}[section]
\theoremstyle{definition} 
\newtheorem{Definition}{Definition}[section]
\newtheorem{Remark}{Remark}[section]
\theoremstyle{plain} 

\newtheorem{Lemma}[Theorem]{Lemma}
\newtheorem{Proposition}[Theorem]{Proposition}

\newtheorem{Corollary}[Theorem]{Corollary}

\def\la{\longrightarrow}
\def\rh{\hookrightarrow}

\DeclareMathOperator*{\colim}{colim}
\DeclareMathOperator*{\hocolim}{hocolim}
\DeclareMathAlphabet{\mathpzc}{OT1}{pzc}{m}{it}
\newcommand{\mf}{{\mathrm{mf}}}

\newcommand{\gmf}{{\mathrm{gmf}}}

\newcommand{\Sets}{{\mathcal S}{\mathrm{ets}}}

\newcommand{\p}{\partial}

\newcommand{\Cob}{{\mathcal C}\!{\mathpzc{o}}\!{\mathpzc{b}}}

\newcommand{\Th}{{\mathsf T}\!{\mathsf h}}

\newcommand{\R}{{\mathbf R}}

\newcommand{\Hom}{{\mathrm{Hom}}}

\newenvironment{Proof}{\par\noindent{\bf Proof}. }{\hfill $\square$\par}

\begin{document}

\title [The moduli space of generalized Morse functions] 
       {The moduli space of generalized Morse functions}

\author{Boris Botvinnik}
\thanks{Boris Botvinnik is partially
  supported by SFB-748, M\"unster, Germany}

\address{Department of Mathematics \\ University of Oregon, Eugene,
  OR, 97403 USA } 

\email{botvinn@math.uoregon.edu}

\author{Ib Madsen} \address{Institute for the Mathematical Sciences,
  University of Copenhagen, DK-2100 Copenhagen, Denmark} 

\thanks{Ib Madsen is partially supported by ERC advanced grant
  228082-TMSS} 

\email{imadsen@math.ku.dk}

\begin{abstract} {We study the moduli and determine a homotopy type of
    the space of all generalized Morse functions on $d$-manifolds for
    given $d$. This moduli space is closely connected to the moduli
    space of all Morse functions studied in \cite{MW} and
    the classifying space of the corresponding
    cobordism category.}
\end{abstract}
\maketitle
\vspace*{-10mm}

\section{Introduction}
Given a smooth compact manifold $M^d$ and a fixed smooth function
$\varphi : M^d \to \R$, let ${\mathcal G}(M^d,\varphi)$ denote the space of
generalized Morse functions $f: M\to \R$ which agrees with $\varphi$ in a
neighbourhood of the boundary $\p M$. This space (in the Whitney
topology) satisfies an $h$-principle in the sense of Gromov,
\cite{G}. Here precisely we define $h{\mathcal G}(M^d,\varphi)$ to be the
space of sections of the bundle $J^3_{\gmf}(M)$ of generalized Morse
$3$-jets that agrees with $j^3\varphi$ near $\p M$. Taking the $3$-jet
of a generalized Morse function defines a map
\begin{equation}\label{eq1}
j^3 : {\mathcal G}(M^d,\varphi) \longrightarrow h{\mathcal G}(M^d,\varphi) \ .
\end{equation}
This map was first considered by Igusa in \cite{I1}. He proved that
the map $j^3$ in (\ref{eq1}) is $d$-connected and in \cite{I2} he calculated
the ``$d$-homotopy type'' of $h{\mathcal G}(M^d,\varphi)$ by exhibiting a
$d$-connected map
\begin{equation}\label{eq2}
h{\mathcal G}(M^d,\varphi)\longrightarrow \Omega^{\infty}S^{\infty}(BO_+\wedge M^d),
\end{equation}
thus determining the $d$-homotopy type of the space ${\mathcal G}(M^d,\varphi)$.

Eliashberg and Mishachev \cite{EM1,EM2} and Vassiliev \cite{V}
showed that the map in (\ref{eq1}) is actually a homotopy equivalence
rather than just being $d$-connected. This is the starting point for
this paper.

We study the moduli space of all generalized Morse functions on
$d$-manifolds, i.e. the space of ${\mathcal G}(M^d,\varphi)$ as
$(M^d,\varphi)$ varies. There can be several candidates for such a moduli
space. The one we present below is closely connected to the ``moduli
space'' of all Morse functions considered in Section 4 of
\cite{MW}. Indeed, the present note can be viewed as an addition to
\cite{MW}. 

In Section 2 below we give the precise definition of our moduli
space, and in Section 3 we determine its homotopy type, following the
argument from \cite{MW}.
\section{Definitions and results}
\subsection{The moduli space} Let $J^3(\R^d)$ be the space of $3$-jets
of smooth functions on $\R^d$,
$$
p(x)= c +\ell(x)+ q(x) + r(x) \ ,
$$
where $c$ is a constant, $\ell(x)$ is linear, $q(x)$ quadratic and
$r(x)$ cubic,
$$
\ell(x)=\sum_i a_ix_i \ , \ \ 
q(x) =\sum_{ij} a_{ij}x_i x_j \ , \ \
r(x) =\sum_{ijk} a_{ijk}x_i x_j x_k \ ,
$$
with the coefficients $a_{ij}$, $a_{ijk}$ 
symmetric in the indices.

Let $J^3_{\gmf}(\R^d)\subset J^3(\R^d)$ be the subspace of $p\in
J^3(\R^d)$ such that one the following holds:
\begin{enumerate}
\item[(i)] $0\in \R^d$ is not a critical point of $p$
  ($\ell\neq 0$);
\item[(ii)] $0\in \R^d$ is a non-degenerate critical point of
  $p$ ($\ell=0$ and $q$ is non-degenerate);
\item[(iii)] $0\in \R^d$ is a birth-death singularity
  of $p$ ($q : \R^d\to \Hom_{\R}(\R^d,\R)$ has $1$-dimensional kernal
  on which $r(x)$ is non-trivial).
\end{enumerate}
The space $J^3_{\gmf}(\R^d)$ is invariant under the
$O(d)$-action on the space $J^3(\R^d)$. Given a smooth manifold $M^d$
with a metric, let ${\mathcal P}(M^d)\to M$ be the principal
$O(d)$-bundle of orthogonal frames in the tangent bundle $TM^d$. Then
$$
J^3_{\gmf}(TM^d)= {\mathcal P}(M^d)\times_{O(d)} J^3_{\gmf}(\R^d)
$$
is a smooth fiber bundle on $M$, a subbundle of
$$
J^3(TM^d)= {\mathcal P}(M^d)\times_{O(d)} J^3(\R^d) .
$$
\begin{Remark}
Up to change of coordinates, a birth-death singularity is of the form
$$
p(x)= x_1^3 - \sum_{j=2}^{i+1} x_j^2 + \sum_{k=i+2}^{d} x_k^2  \ .
$$
The integer $i$ is the Morse index of the quadratic form $q(x)$. In
general, a quadratic form $q: \R^d\to \R$ induces a canonical decomposition
$$
\R^d=V_-(q)\oplus V_0(q)\oplus V_+(q)
$$
into negative eigenspace, the zero eigenspace and the positive
eigenspace. In the case of generalized Morse jets, $\dim V_0(q)$ is
either $0$ or $1$, and in the latter case the cubic term $r(x)$
restricts non-trivially to $V_0(q)$. The dimension of $V_-(q)$ is the
index of the gmf-jet. \hfill $\Diamond$
\end{Remark}
For a smooth manifold $M^d$, we have the $3$-jet bundle $J^3(M,\R)\to
M$ whose fiber $J^3(M,\R)_x$ is the germ of $3$-jets $(M,x)\to \R$,
and the associated subbundle $J^3_{\gmf}(M,\R)\to M$. A choice of
exponential function induces a fiber bundle isomorphism
\begin{equation}\label{eq4}
J^3_{\gmf}(M,\R) \cong J^3_{\gmf}(TM),
\end{equation}
and $f : M\to \R$ is a generalized Morse function precisely if $j^3(f)\in
\Gamma(J^3_{\gmf}(M,\R) )$, where we use $\Gamma(E)$ to denote the
space of smooth sections of a vector bundle $E$.
\begin{Definition}\label{def1}
  Let $X$ be a $k$-dimensional manifold without boundary. Let
  ${\mathcal J}_d(X)$ be the set of $4$-tuples $(E,\pi,f,j)$ of a
  $(k+d)$-manifold $E$ with maps 
$$
(\pi,f,j): E\longrightarrow X\times \R\times \R^{d-1+\infty}
$$
subject to the conditions
\begin{enumerate}
\item[(i)] $(\pi,f): E \longrightarrow X\times \R$ is
  a proper map;
\item[(ii)] $(f,j): E\longrightarrow \R\times
  \R^{d-1+\infty}$ is an embedding;
\item[(iii)] $\pi : E\longrightarrow X$ is a
  submersion;
\item[(iv)] for any $x\in X$, the restriction
  $f_x=f|_{E_x}: E_x \longrightarrow \R$ to each fiber
    $E_x=\pi^{-1}(x)$ is a generalized Morse function.
\end{enumerate}
\end{Definition}
In (ii) above $\R^{d-1+\infty}$ is the union or colimit of
$\R^{d-1+N}$ as $N\to \infty$ and (ii) means that $(f,j)$ embeds $E$
into $\R\times \R^{d-1+N}$ for sufficiently large $N$.  The definition
above is the obvious analogue of Definition 2.7 of \cite{MW}.

A smooth map $\phi : Y\longrightarrow X$ induces a pull-back
$$
\begin{array}{c}
  \phi^* : {\mathcal J}_d(X)\la {\mathcal J}_d(Y), \ \ \ 
  \phi^*: (E,\pi,f,j)\mapsto (\phi^*E,\phi^*\pi,\phi^*f,\phi^*j)
  \ \ \ \mbox{where}
  \\
  \\
  \phi^*E = \{ \ (y,z)\in Y\times \R^{d+\infty} \ | 
  \ (\phi(y),z)\in E\subset X\times \R^{d+\infty}\ \}
\end{array}
$$
and the maps $\phi^*\pi$, $\phi^*f$ and $\phi^*j$ are given by
corresponding projections from 
$$
\phi^*E \subset Y\times \R\times \R^{d-1+\infty}
$$ 
on the factors $Y$, $\R$ and $\R^{d-1+\infty}$,
respectively.  In particular, we have
$(\psi\circ\phi)^*=\phi^*\circ\psi^*$ (rather than just being
naturally equivalent), so the correspondence $X\mapsto {\mathcal
  J}_d(X)$ is a set-valued sheaf ${\mathcal J}_d$ on the category
${\mathfrak X}$ of smooth manifolds and smooth maps.

A set-valued sheaf on ${\mathfrak X}$ gives rise to a simplicial set
$N_{\bullet}{\mathcal J}_d$ with
$$
\begin{array}{c}
N_{k}{\mathcal J}_d = {\mathcal J}_d(\Delta^k_e), \ \ \
\Delta^k_e = \left\{ \ (x_0,\ldots,x_k)\in \R^{k+1} \ | \ \sum x_i =1 \ \right\}.
\end{array}
$$
The geometric realization of $N_{\bullet}{\mathcal J}_d$ will be denoted by
$|{\mathcal J}_d|$, and we make the following definition.
\begin{Definition}\label{def2}
  The moduli space of generalized Morse functions of $d$ variables is
  the loop space $\Omega|{\mathcal J}_d|$.
\end{Definition}
\begin{Remark}
  If in Definition \ref{def1} we drop the assumption that $f: E\la \R$
  is a generalized Morse function, then ${\mathcal J}_d$ reduces to
  the sheaf $D_d=D_d(-,\infty)$ of \cite[Definition 3.3]{GMTW}
  associated to the space of embedded $d$-manifolds, and
  $\Omega|D_d|\cong \Omega^{\infty} MT(d)$ by Theorem 3.4 of
  \cite{GMTW}. 
  \hfill $\Diamond$
\end{Remark}

Associated with the set valued sheaf ${\mathcal J}_d(X)$, $X\in
{\mathfrak X}$, we have a sheaf ${\mathcal J}_d^{{\mathcal A}}(X)$ of
partially ordered sets, i.e. a category valued sheaf,
cf. \cite[Section 4.2]{MW}. For connected $X$, an object of ${\mathcal
  J}_d^{{\mathcal A}}(X)$ consists of an element $(E,\pi,f,j)\in
{\mathcal J}_d(X)$ together with an interval $A=[a_0,a_1]\subset \R$
subject to the condition that $f: E\la \R$ be fiberwise transverse to
$\p A$ (i.e. $\{a_0,a_1\}$ are regular values for each 
$f_x: E_x \la \R$, $x\in X$). The partial ordering is given by
$$
\begin{array}{c}
(E,\pi,f,j;A) \leq (E',\pi',f',j';A') \ \ \mbox{if}
\\
\\
(E,\pi,f,j) = (E',\pi',f',j')  \ \ \mbox{and} \ \ A\subset A' \ .
\end{array}
$$
If $X$ is not connected, $\displaystyle {\mathcal J}_d^{{\mathcal
    A}}(X)$ is the product of ${\mathcal J}_d^{{\mathcal A}}(X_j)$
over the connected components $X_j$.

We notice that each element $(E,\pi,j,f;A)$ restricts to a family of
generalized Morse functions on a compact manifolds
$$
(\pi,f,j) : f^{-1}(A)\rh X\times [a_0,a_1] \times \R^{d-1+\infty},
$$
where $A=[a_0,a_1]$. On the other hand, given such a family
$$
(\pi,f,j) : E(A)\rh X\times [a_0,a_1] \times \R^{d-1+\infty},
$$
we can extend it to an element $(\hat E(A), \hat \pi,\hat f,\hat j)$
by adding long collars:
$$
\hat E(A) = (-\infty, a_0]\times f^{-1}(a_0)\cup E(A)\cup 
[a_1,\infty)\times f^{-1}(a_1).
$$
If $(E(A)\pi,f,j)$ is a restriction of $(E,\pi,f,j)\in {\mathcal
  J}_d(X)$, then $(\hat E(A), \hat \pi,\hat f,\hat j; A)$ is concordant to 
$(E,\pi,f,j)$ by \cite[Lemma 2.19]{MW}. 

The forgetful map ${\mathcal J}_d^{{\mathcal A}}(X)\la {\mathcal
  J}_d(X) $ is a map of category valued sheaves when we give
${\mathcal J}_d(X)$ the trivial category structure (with only identity
morphisms). It induces a map
$$
|{\mathcal J}_d^{{\mathcal A}}|\la |{\mathcal J}_d|
$$
of topological categories, and hence a map of their classifying spaces:
$$
B|{\mathcal J}_d^{{\mathcal A}}|\la B|{\mathcal J}_d|=|{\mathcal J}_d|,
$$
where $B|{\mathcal J}|= |N_{\bullet}{\mathcal J}|$.
\begin{Theorem}\label{thm2.5}
  The map $B|{\mathcal J}_d^{{\mathcal A}}|\la |{\mathcal J}_d|$ is a
  weak homotopy equivalence. 
\end{Theorem}
\begin{Proof}
  This follows from \cite[Theorem 4.2]{MW}, which identifies $|B
  {\mathcal J}_d^{{\mathcal A}}|$ with $|\beta {\mathcal
    J}_d^{{\mathcal A}}|$, where $\beta {\mathcal J}_d^{{\mathcal A}}$
  is a set-valued sheaf of \cite[Definition 4.1]{MW}, together with
  the analogue of \cite[Proposition 4.10]{MW}: the map
$$
|\beta {\mathcal J}_d^{{\mathcal A}}| \la |{\mathcal J}_d^{{\mathcal
    A}}|
$$
is a weak homotopy equivalence. Indeed, the proof of Proposition 4.10,
which treats the case where $f: E \la \R$ is a fiberwise Morse
function (in a neighbourhood of $f^{-1}(0)$) carries over word by word
to the situation of generalized Morse functions.
\end{Proof}
\subsection{The $h$-principle} For a submersion $\pi: E \la X$,
$T^{\pi}E$ denotes the tangent bundle along the
fibers. We can form the bundle
$$
J^3_{\gmf}(T^{\pi}E) \la E
$$
of gmf-jets. Sections of this bundle will be denoted by $\hat f$,
$\hat g,\ldots$ etc. For given $z\in E$, the restriction $\hat f(z) :
T(E_{\pi(z)})\la \R$ is a gmf-jet; its constant term will be denoted
$f(z)$.
\begin{Definition}\label{def2.6}
For a smooth manifold $X$, let $h{\mathcal J}_d(X)$ consists of maps
$$
(\pi,f,j): E\longrightarrow X\times \R\times \R^{d-1+\infty}
$$
satisfying (i), (ii), and (iii) of Definition \ref{def1} together with
a jet $\hat f\in \Gamma(J^3_{\gmf}(T^{\pi}E))$ having constant term
$f$. \hfill $\Box$
\end{Definition}
A metric on $T^{\pi}E$ and an associated exponential
map induces an isomorphism
\begin{equation}\label{eq5}
J^3_{\pi}(E,\R) \stackrel{\cong}{\la} J^3(T^{\pi}E)
\end{equation}
that sends gmf-jets to gmf-jets. Here $J^3_{\pi}(E,\R) \la E$ is the
fiberwise $3$-jet bundle. Differentiation in the fiber direction only,
defines a map
$$
j^3_{\pi}: C^{\infty}(E,\R) \la J^3_{\pi}(E,\R)
$$
which sends fiberwise generalized Morse functions into gmf-jets. This
induces a map of sheaves
$$
j^3_{\pi} : {\mathcal J}_d(X)\la h{\mathcal J}_d(X),
$$
and hence a map of their nerves
$$
j^3_{\pi} :  |{\mathcal J}_d|\la |h{\mathcal J}_d|.
$$
Using the sheaves ${\mathcal J}_d^{{\mathcal A}}$ and the associated
$h{\mathcal J}_d^{{\mathcal A}}$, together with the $h$-principle of
\cite{EM2,V}, the argument of \cite[Proposition 4.17]{MW} proves that
the induced map
$$
|\beta {\mathcal J}_d^{{\mathcal A}}|\la |\beta h{\mathcal
  J}_d^{{\mathcal A}}|
$$
is a weak homotopy equivalence. Finally, since the forgetful maps
$$
|\beta {\mathcal
  J}_d^{{\mathcal A}}| \la |{\mathcal
  J}_d^{{\mathcal A}}|, \ \ \ \ 
|\beta h{\mathcal
  J}_d^{{\mathcal A}}| \la |h{\mathcal
  J}_d^{{\mathcal A}}|
$$
are weak homotopy equivalences by \cite[Proposition 4.10]{MW}, we get
\begin{Theorem}\label{thm2-7}
  The map $j^3_{\pi} : |{\mathcal J}_d|\la |h{\mathcal J}_d|$ is weak
  homotopy equivalence. \hfill $\Box$
\end{Theorem}
\section{Homotopy type of the moduli space}
\subsection{The space $|h{\mathcal J}_d|$} 
For any set-valued sheaf ${\mathcal F} : {\mathfrak X} \la \Sets$, let
${\mathcal F}[X]$ denote the set of concordance classes: $s_0,s_1\in
{\mathcal F}(X)$ are \emph{concordant} ($s_0\sim s_1$) if there exists
an element $s\in {\mathcal F}(X\times \R)$ such that $pr_X^*(s_0)$,
$pr_X : X\times \R \la X$, agrees with $s$ on an open neighborhood of
$X\times (-\infty,0]$ and $pr_X^*(s_1)$ agrees with $s$ on an open
neighborhood of $X\times [1,\infty)$. The relation to the space
$|{\mathcal F}|$ is given by
\begin{equation}\label{eq7}
[X, |{\mathcal F}|]\cong {\mathcal F}[X].
\end{equation}
Let $G(d,n)$ denote the Grassmannian of $d$-planes in $\R^{d+n}$, and
$G^{\gmf}(d,n)$ the space of pairs $(V,f)$ with $V\in G(d,n)$ and $f:
V\to \R$ a generalized Morse function with $f(0)=0$. The space
$$
U^{\perp}_{d,n}= \{ \ (v,V) \in \R^{d+n}\times G(d,n) \ | \ v\perp V \ \}
$$
is an $n$-dimensional vector bundle on $G(d,n)$. Let $V^{\perp}_{d,n}$
be its pull-back along the forgetful map $G^{\gmf}(d,n)\la G(d,n)$:
$$
\begin{diagram}
\setlength{\dgARROWLENGTH}{1.2em}
       \node{V^{\perp}_{d,n}}
           \arrow[2]{e}
           \arrow{s}
       \node[2]{U^{\perp}_{d,n}}
           \arrow{s}
\\
       \node{G^{\gmf}(d,n)}
           \arrow[2]{e}
       \node[2]{G(d,n)}
\end{diagram}
$$
Similarly, we have canonical $d$-dimensional vector bundles $U_{d,n}$
and $V_{d,n}$ on $G(d,n)$ and $G^{\gmf}(d,n)$ respectively. The Thom
spaces of the bundles $U^{\perp}_{d,n}$ and $V^{\perp}_{d,n}$ give
rise to spectra $MT(d)$ and $MT^{\gmf}(d)$ which in degrees $(d+n)$
are
\begin{equation}\label{eq8b}
MT(d)_{d+n} = \Th(U^{\perp}_{d,n}), \ \ \ \ 
MT^{\gmf}(d)_{d+n}  = \Th(V^{\perp}_{d,n}).
\end{equation}
The infinite loop space of the spectrum $MT^{\gmf}(d)$ is defined to
be 
\begin{equation}\label{eq7}
\Omega^{\infty}MT^{\gmf}(d)=\colim_{n} \Omega^{d+n} \Th(V^{\perp}_{d,n})
\end{equation}
\begin{Theorem}\label{thm3-1}
There is a weak homotopy equivalence
\begin{equation}\label{eq7a}
\Omega
|h{\mathcal J}_d| \simeq \Omega^{\infty}MT^{\gmf}(d).
\end{equation}
\end{Theorem}
\begin{Proof}
  This is completely similar to the proof of \cite[Theorem 3.5]{MW} for
  the case of Morse functions using only transversality and the
  submersion theorem, \cite{Phillips}.
\end{Proof}
We have left to examine the right-hand side of the equivalence
(\ref{eq7a}). The results are similar in spirit to ones in
\cite[Section 3.1]{MW}.

Let $\Sigma^{\gmf}(d,n)\subset G^{\gmf}(d,n)$ be the singular set of
pairs $(V,f)$ with $f: V\la \R$ having vanishing linear term,
i.e. $Df(0)=0$.  Consider the nonsingular set $G^{\gmf}(d,n)\setminus
\Sigma^{\gmf}(d,n)$ given as set of pairs $(V,f)\in G^{\gmf}(d,n)$
with $f=\ell+q+r$ and $\ell\neq 0$. We notice that the space
$G^{\gmf}(d,n)\setminus \Sigma^{\gmf}(d,n)$ retracts to the space
$$
\hat G^{\gmf}(d,n) = \{ \ (V,f) \in G^{\gmf}(d,n) \ | \ f=\ell+q+r, \ 
\mbox{with} \ |\ell|= 1, \ q=0, \ r=0 \},
$$
where $|\ell|$ is a norm of the linear part $\ell: V\la \R$. 
\begin{Lemma}\label{lem-new1}
There is a homeomorphism 
$$
\hat G^{\gmf}(d,n) \cong 
O(d+n)/\left(O(d-1)\times O(n)\right) \ .
$$
\end{Lemma}
\begin{Proof}
  For a pair $(V,f)\in \hat G^{\gmf}(d,n)$ we have $V\in G(d,n)$ and
  $f=\ell$ with $|\ell|=1$. We may think of $\ell$ as a linear
  projection on the first coordinate, which is the same as to say that
  the space $V$ contains a subspace 
\begin{equation}\label{eq-new1}
\R\times \{0\}\times \{0\}\subset \R\times \R^{d-1}\times \R^n
\end{equation}
with $\ell$ being a projection on it. This identifies $\hat
G^{\gmf}(d,n)$ with the homogeneous space $ O(d+n)/\left(O(d-1)\times
  O(n)\right)$.
\end{Proof}
Since $G(d-1,n)=O(d-1+n)/(O(d-1)\times O(n))$, we
observe that the map
$$
i_{d,n}: G(d-1,n) \longrightarrow \hat G^{\gmf}(d,n)
$$
is $(d+n-2)$-connected and that
$i_{d,n}^*(V_{d,n}^{\perp}|_{\hat
      G^{\gmf}(d,n)})\cong U_{d-1,n}^{\perp}$.

On the other hand, $\Sigma^{\gmf}(d,n)\subset G^{\gmf}(d,n)$ has normal
bundle $V_{d,n}^*\cong V_{d,n}$ and the inclusion
$$
(D(V_{d,n}),S(V_{d,n})) \la (G^{\gmf}(d,n), G^{\gmf}(d,n)\setminus
\Sigma^{\gmf}(d,n))
$$
is an excision map. This leads to the cofibration
\begin{equation}\label{eq9}
  \Th(j^*V^{\perp}_{d,n})\la \Th(V^{\perp}_{d,n})\la 
\Th((V^{\perp}_{d,n}\oplus V_{d,n})|_{\Sigma^{\gmf}(d,n)}),
\end{equation}
where $j$ is the inclusion 
$$
j : G^{\gmf}(d,n)\setminus \Sigma^{\gmf}(d,n) \la G^{\gmf}(d,n) \ .
$$
By the above, there is $(2n+d-2)$-connected map
$$
i_{d,n}: \Th(U^{\perp}_{d-1,n}) \la
\Th(V^{\perp}_{d,n}) .
$$
With the notation of (\ref{eq8b}), we get from (\ref{eq9}) a
cofibration of spectra
\begin{equation}\label{eq10}
\Sigma^{-1}MT(d-1) \la MT^{\gmf}(d)\la \Sigma^{\infty}(\Sigma^{\gmf}(d,\infty)_+) 
\end{equation}
and the corresponding homotopy fibration sequence of infinite loop
spaces
$$
\Omega^{\infty}
\Sigma^{-1}MT(d-1) \la \Omega^{\infty}MT^{\gmf}(d)\la 
\Omega^{\infty}\Sigma^{\infty}(\Sigma^{\gmf}(d,\infty)_+) 
$$
\begin{Remark} The main theorem of \cite{GMTW}
asserts a homotopy equivalence
$$
\Omega^{\infty} MT(d) \simeq  \Omega B\Cob_d  ,
$$
where $\Cob_d$ is the cobordism category of embedded manifolds: the
objects are $(M^{d-1},a)$ of a closed $(d-1)$-submanifold of
$\{a\}\times \R^{\infty+d-1}$ and the morphisms are embedded
cobordisms $W^d\subset [a_0,a_1]\times \R^{\infty+d-1}$ transversal at
$\{a_i\}\times \R^{\infty+d-1}$. In particular, we have weak
homotopy equivalence
\begin{equation}\label{eq11}
\Omega^{\infty}
\Sigma^{-1}MT(d-1) \simeq \Omega^2 B \Cob_{d-1}.
\end{equation}
\end{Remark}
\subsection{The singularity space}
In \cite{I2}, Igusa analyzed the singularity space
$\Sigma^{\gmf}(d)\subset J^3_{\gmf}(\R^d)$ by
  decomposing it with respect to the Morse index. The result, stated
  in \cite[Proposition 3.4]{I2}, is as follows.  Consider the
  homogeneous spaces 
$$
\begin{array}{lcl}
  X^1(i) & = & \displaystyle
  O(d) \left/\mbox{\begin{picture}(1,11)\end{picture}}O(i)\times O(1)\times O(d-i-1) , \right.
  \\
  \\
  X(i) & = & \displaystyle
  O(d) \left/\mbox{\begin{picture}(1,11)\end{picture}} 
O(i)\times O(d-i) .\right.
\end{array}
$$
and note that there are quotient maps
$$
f_i : X^1(i) \to X(i), \ \ \ \ g_i : X^1(i) \to X(i+1),
$$
upon embedding $O(i)\times O(1)\times O(d-i-1)$ in $O(i)\times O(d-i)$
and in $O(i+1)\times O(d-i-1)$, respectively. These maps fit into the 
diagram 
$$
{\mathcal D}(d) =\left(\!\!
\begin{diagram}
\setlength{\dgARROWLENGTH}{-1.9em}
\node[2]{\!\!X^1(0)\!\!}
     \arrow{sse,t}{g_0}
     \arrow{ssw,t}{f_0}
\node[2]{\!\!X^1(1)\!\!}
     \arrow{sse,t}{g_1}
     \arrow{ssw,t}{f_1}
\node[2]{\!\!\cdots\!\!}
     \arrow{sse,t}{g_{d-2}}
     \arrow{ssw,t}{f_2}
\node[2]{\!\!X^1(d-1)\!\!}
     \arrow{sse,t}{g_{d-1}}
     \arrow{ssw,t}{f_{d-1}}
\\
\\
\node{X(0)\!\!\!\!}
\node[2]{\!\!X(1)\!\!}
\node[2]{\!\!X(2)\!\!}
\node[2]{\!\!X(d-1)\!\!}
\node[2]{\!\!X(d)}
\end{diagram} \!\!
\right)
$$
and \cite[Proposition 3.4]{I2} states that the homotopy colimit of the
diagram ${\mathcal D}(d)$ is homotopy equivalent to $\Sigma^{\gmf}(d)$
\begin{equation}\label{eq12}
\Sigma^{\gmf}(d) \simeq \hocolim {\mathcal D}(d).
\end{equation}
It is easy to see that there are homeomorphisms 
$$
\begin{array}{lcl}
  G^{\gmf}(d,n) &=& \left(
    O(d+n) \left/\mbox{\begin{picture}(1,11)\end{picture}} O(n)\right. 
    \right)\times_{O(d)} J^3_{\gmf}(\R^d),
  \\
  \\
  \Sigma^{\gmf}(d,n) &=& \left( O(d+n)
   \left/\mbox{\begin{picture}(1,11)\end{picture}}  O(n)\right. 
\right)\times_{O(d)} \Sigma^{\gmf}(d),
\end{array}
$$
and (\ref{eq12}) implies that $\Sigma^{\gmf}(d,n)$ is
homotopy equivalent to the homotopy colimit of the diagram
$$
{\mathcal D} (d,n)= \left( O(d+n)
  \left/\mbox{\begin{picture}(1,11)\end{picture}} O(n)\right.
\right)\times_{O(d)} {\mathcal D} (d) \ .
$$
For $n\to\infty$, the Stiefel manifold $O(n+d)/O(n)$ becomes
contractible, and ${\mathcal D}(d,\infty)$ is the diagram
\begin{equation}\label{eq13}
\begin{diagram}
\setlength{\dgARROWLENGTH}{-1.9em}
\node[2]{\!\!Y^1(0)\!\!}
     \arrow{sse,t}{\bar g_0}
     \arrow{ssw,t}{\bar f_0}
\node[2]{\!\!Y^1(1)\!\!}
     \arrow{sse,t}{\bar g_1}
     \arrow{ssw,t}{f_1}
\node[2]{\!\!\cdots\!\!}
     \arrow{sse,t}{\bar g_{d-2}}
     \arrow{ssw,t}{\bar f_2}
\node[2]{\!\!Y^1(d-1)\!\!}
     \arrow{sse,t}{\bar g_{d-1}}
     \arrow{ssw,t}{\bar f_{d-1}}
\\
\\
\node{Y(0)\!\!\!\!}
\node[2]{\!\!Y(1)\!\!}
\node[2]{\!\!Y(2)\!\!}
\node[2]{\!\!Y(d-1)\!\!}
\node[2]{\!\!Y(d)}
\end{diagram} \!\!
\end{equation}
with
$$
\begin{array}{lcl}
Y^1(i) & = & BO(i)\times BO(1)\times BO(d-i-1),
\\
\\
Y(i) & = & BO(i)\times BO(d-i),
\end{array}
$$
and $\bar f_i$ and $\bar g_i$ the obvious maps. So
$\Sigma^{\gmf}(d,\infty)$ is the homotopy colimit of (\ref{eq13}).

We want to compare this to the singular set $\Sigma^{\mf}(d,\infty)$
which appears when one considers the moduli space of Morse functions
rather than generalized Morse functions was calculated in \cite[Lemma
3.1]{MW}:
$$
\Sigma^{\mf}(d,n) \cong \prod_{i=0}^d \left[ \left( O(d+n)
    \left/\mbox{\begin{picture}(1,11)\end{picture}} O(n)\right.
  \right) \times_{O(d)} \left(O(d)
    \left/\mbox{\begin{picture}(1,11)\end{picture}} O(i)\times
      O(d-i)\right. \right)\right]
$$
so that 
$$
\Sigma^{\mf}(d,\infty) \cong \prod_{i=0}^d BO(i)\times BO(d-i) = 
\prod_{i=0}^d Y(i).
$$
The cofiber of the map $\Sigma^{\mf}(d,\infty)\la
\Sigma^{\gmf}(d,\infty)$ is by (\ref{eq13}) equal to the homotopy
colimit of the diagram:
$$
\begin{diagram}
\setlength{\dgARROWLENGTH}{0.6em}
\node[2]{\!\!Y^1(0)\!\!}
     \arrow{sse}
     \arrow{ssw}
\node[2]{\!\!Y^1(1)\!\!}
     \arrow{sse}
     \arrow{ssw}
\node[2]{\!\!\cdots\!\!}
     \arrow{sse}
     \arrow{ssw}
\node[2]{\!\!Y^1(d-1)\!\!}
     \arrow{sse}
     \arrow{ssw}
\\
\\
\node{*}
\node[2]{*}
\node[2]{*}
\node[2]{*}
\node[2]{*}
\end{diagram} \!\!
$$
But this homotopy colimit is easy calculated to be
$$
\bigvee_{i=0}^{d-1} S^1\wedge Y^1(i)_+\simeq
\bigvee_{i=0}^{d-1} S^1\wedge (BO(i)\times BO(d-i-1))_+ \ .
$$
We get a cofibration of suspension spectra:
$$
\Sigma^{\infty}(\Sigma^{\mf}(d,\infty)_+) \la
\Sigma^{\infty}(\Sigma^{\gmf}(d,\infty)_+) \la \bigvee_{i=0}^{d-1}
\Sigma^{\infty}(S^1\wedge
(BO(i)\times BO(d-i-1)_+)).
$$
Taking the associated infinite loop spaces we get
\begin{Proposition}\label{prop3-2}
There is a homotopy fibration:
$$
\prod_{i=0}^{d-1} \Omega^{\infty}\Sigma^{\infty}(BO(i)\!\times 
\!BO(d\!-\!i\!-\!1)_+)
\to
\prod_{i=0}^{d} \Omega^{\infty}\Sigma^{\infty}
(\Sigma^{\mf}(d,\infty)_+)
\to 
\Omega^{\infty}\Sigma^{\infty}(\Sigma^{\gmf}(d,\infty)_+). \ 
\Box \!\!\!\!
$$
\end{Proposition}
The constant map ${\mathcal D}(d)\to *$ into the constant diagram
induces the map
$$
{\mathcal D}(d,n) \to \left( O(d+n)
    \left/\mbox{\begin{picture}(1,11)\end{picture}} O(n)\right.
  \right) \times_{O(d)} * \ , 
$$
where the target space is homotopy equivalent to $G(d,n)$.  For
$n\to\infty$, this induces the fiber bundle
\begin{equation}\label{bo-d}
p: \Sigma^{\gmf}(d,\infty) \to BO(d)
\end{equation}
with the fiber $\Sigma^{\gmf}(d)$. We obtain the commutative diagram
of cofibrations:
$$
\begin{diagram}
\setlength{\dgARROWLENGTH}{1.2em}
\node{\Sigma^{-1}MT(d-1)}
     \arrow[2]{e}
     \arrow{s,l}{Id}
\node[2]{MT^{\gmf}(d)}
     \arrow[2]{e}
     \arrow{s,l}{F}
\node[2]{\Sigma^{\infty}(\Sigma^{\gmf}(d,\infty)_+)}
     \arrow{s,l}{\Sigma^{\infty}p}
\\
\node{\Sigma^{-1}MT(d-1)}
     \arrow[2]{e}
\node[2]{MT(d)}
     \arrow[2]{e}
\node[2]{\Sigma^{\infty}(BO(d)_+)}
\end{diagram} 
$$
and a corresponding diagram of homotopy fibrations:
\begin{equation}\label{gmf-new}
\begin{diagram}
\setlength{\dgARROWLENGTH}{1.2em}
\node{\Omega^{\infty}\Sigma^{-1}MT(d-1)}
     \arrow[2]{e}
     \arrow{s,l}{Id}
\node[2]{\Omega^{\infty}MT^{\gmf}(d)}
     \arrow[2]{e}
     \arrow{s,l}{\Omega^{\infty} F}
\node[2]{\Omega^{\infty}\Sigma^{\infty}(\Sigma^{\gmf}(d,\infty)_+)}
     \arrow{s,l}{\Omega^{\infty}\Sigma^{\infty}p}
\\
\node{\Omega^{\infty}\Sigma^{-1}MT(d-1)}
     \arrow[2]{e}
\node[2]{\Omega^{\infty}MT(d)}
     \arrow[2]{e}
\node[2]{\Omega^{\infty}\Sigma^{\infty}(BO(d)_+)}
\end{diagram} 
\end{equation}
Since $\Sigma^{\infty}(\Sigma^{\gmf}(d,\infty)_+)$ and
$\Sigma^{\infty}(BO(d)_+)$ are $(-1)$-connected, we obtain that the
forgetful map $F: MT^{\gmf}(d)\to MT(d)$ induces isomorphism
$$
\pi_{-i}MT^{\gmf}(d)\cong
\pi_{-i}MT(d), \ \ \ i\geq 0 .
$$
We consider the forgetful map $ \theta^{\gmf}:
  G^{\gmf}(d,\infty) \longrightarrow G(d,\infty) $ as a
  \emph{structure on $d$-dimensional bundles}. Then we denote by
  $\Cob^{\gmf}_d$ the category $\Cob^{\theta^{\gmf}}_d$ (see (5.3) and
  (5.4) of \cite{GMTW}) of manifolds (objects) and cobordisms
  (morphisms) equipped with a tangential $\theta^{\gmf}$-structure.
  Then the main theorem of \cite{GMTW} gives the following result:
  \begin{Corollary}\label{cor1-new}
There is weak homotopy equivalence
$$
B\Cob_d^{\gmf}\cong \Omega^{\infty-1}MT^{\gmf}(d),
$$
and the forgetting map $B\Cob_d^{\gmf}\to B\Cob_d$ induces isomorphism
$$
\Omega^{\gmf}_d=\pi_0 B\Cob_d^{\gmf}\simeq \pi_0
B\Cob_d=\Omega_d,
$$
where $\Omega^{\gmf}_d$ and $\Omega_d$ are corresponding cobordism groups.
\end{Corollary}
\subsection{Remarks on the moduli space of Morse functions} The paper
\cite{MW} studied the moduli space of fiberwise Morse functions. The
fibers are the space of functions which locally has $2$-jets of the
form $f: \R^{d}\to \R$, $f=f(0)+\ell(x)+q(x)$ subject to the
conditions:
\begin{enumerate}
\item[(i)] $f(0)\neq 0$ or
\item[(ii)] $f(0)=0$ and $\ell(x)\neq 0$ or
\item[(iii)] $f(0)=0$, $\ell(x)=0$ and $q(x)$ is non-singular
  quadratic form.
\end{enumerate}
The associated sheaf ${\mathcal J}^{\mf}_d(X)$, denoted by ${\mathcal
  W}(X)$ in \cite{MW}, consists of maps
$$
(\pi,f,j): E\la X\times \R\times \R^{d-1+\infty} \ \ \ \ \mbox{with}
$$
\begin{enumerate}
\item[(a)] $(\pi,f)$ is proper map,
\item[(b)] $(f,j)$ is an embedding,
\item[(c)] $\pi : E\la X$ is a submersion of relative
  dimension $d$,
\item[(d)] for $x\in X$, $f_x : E_x \to \R$ is ``Morse'', i.e. its
  2-jet satisfies the conditions (i), (ii) and (iii).
\end{enumerate}
The space $|{\mathcal J}^{\mf}_d|=|{\mathcal W}|$ was determined up to
homotopy in Theorems 1.2 and 3.5 of \cite{MW}. We recall the results.
Let $G^{\mf}(d,n)$ be the space of pairs
$$
(V,f) \in G(d,n)\times J^2(V)
$$
with $f$ satisfying the above conditions (i), (ii) and (iii) and
$f(0)=0$. Let $\hat U^{\perp}_{d,n}$ be the canonical $n$-dimensional
bundle over $G^{\mf}(d,n)$ and $MT^{\mf}(d)$ be the spectrum with
$$
MT^{\mf}(d)_{d+n}= \Th(\hat U^{\perp}_{d,n}) \ .
$$
\begin{Theorem}\label{thm3-3}{\rm ([MW])} There is a homotopy
    equivalence
$$
\Omega |{\mathcal J}^{\mf}_d|\cong \Omega^{\infty} MT^{\mf}(d):=
\colim_{n}\Omega^{d+n}\Th(\hat U^{\perp}_{d,n}) \ .
$$
\end{Theorem}
Analogous to (\ref{eq10}), there is the cofibration of spectra
\begin{equation}\label{eq15}
\Sigma^{-1}MT(d-1) \la 
MT^{\mf}(d)\la \Sigma^{\infty}(\Sigma^{\mf}(d,\infty)_+) 
\end{equation}
The inclusion ${\mathcal J}^{\mf}_d(X)\la {\mathcal J}^{\gmf}_d(X)$
induces a map 
$$
\Omega |{\mathcal J}^{\mf}_d|\la \Omega |{\mathcal J}^{\gmf}_d|
$$
of moduli spaces which can be examined upon comparing the (\ref{eq10})
and (\ref{eq15}) We have the homotopy commutative diagram of homotopy
fibrations
\begin{equation}\label{eq16}
\begin{diagram}
\setlength{\dgARROWLENGTH}{1.2em}
\node{\Omega^{\infty}\Sigma^{-1}MT(d-1)}
      \arrow{s,r}{\cong}
      \arrow[2]{e}
\node[2]{\Omega^{\infty} MT^{\mf}(d)}
      \arrow{s}
      \arrow[2]{e}
\node[2]{\Omega^{\infty} \Sigma^{\infty} (\Sigma^{\mf}(d,n)_+)}
      \arrow{s}
\\
\node{\Omega^{\infty}\Sigma^{-1}MT(d-1)}
      \arrow[2]{e}
\node[2]{\Omega^{\infty} MT^{\gmf}(d)}
      \arrow[2]{e}
\node[2]{\Omega^{\infty} \Sigma^{\infty} (\Sigma^{\gmf}(d,n)_+)}
\end{diagram}
\end{equation}
The middle vertical row can be identified with the map
$$
\Omega |{\mathcal J}^{\mf}_d| \la \Omega |{\mathcal J}^{\gmf}_d| 
$$
and the right-hand vertical row corresponds to the right-hand arrow of
Proposition \ref{prop3-2}.  
This gives 
\begin{Corollary}\label{cor3-4}
There is a homotopy fibration
$$
\ \ \ \ \ \ \ \ \ \ \ \  \ \ \ \ \ \ \ \prod_{i=0}^{d-1}
\Omega^{\infty}\Sigma^{-1}(BO(i)\times BO(d-i-1))_+ \la \Omega
|{\mathcal J}^{\mf}_d| \la \Omega |{\mathcal J}^{\gmf}_d| \ . \ \ \ \
\ \ \ \ \ \ \ \ \  \ \ \ \ \ \ \ \ \Box
$$
\end{Corollary}
\subsection{Generalization to tangential structures} Let $\theta :
B\to BO(d)$ be a Serre fibration thought as a
\emph{structure on $d$-dimensional vector bundles}: If $f: X \to
BO(d)$ is a map classifying a vector bundle over $X$, then a map
$\ell: X \to B$ such that $f = \theta \circ \ell$. For a given
$n$, we define the space $G^{\theta,\gmf}(d,n)$ as the pull-back:
\begin{equation}\label{theta}
\begin{diagram}
\setlength{\dgARROWLENGTH}{1.2em}
\node{G^{\theta,\gmf}(d,n)}
     \arrow[2]{e}
     \arrow{s,r}{\theta_{d,n}}
\node[2]{B}
     \arrow{s,r}{\theta}
\\
\node{G^{\gmf}(d,n)}
     \arrow[2]{e,t}{i_{n}}
\node[2]{BO(d)}
\end{diagram} 
\end{equation}
where $i_n$ is the composition of the forgetful map $G^{\gmf}(d,n)\to
G(d,n)$ and the canonical embedding $i^o_n : G(d,n)\rh
G(d,\infty)=BO(d)$.  To define a corresponding sheaf ${\mathcal
  J}^{\theta}_d$ of generalized Morse function on manifolds with
tangential structure $\theta$, we use Definition \ref{def1} but adding
the requirement that the manifold $E$ and corresponding fibers
$E_x=\pi^{-1}(x)$ are equipped with the compatible tangential
structures. Similarly the sheaf $h{\mathcal J}^{\theta}_d$ is
well-defined and there is the corresponding map $j^3_{\pi}(\theta) :
|{\mathcal J}_d^{\theta}|\la |h{\mathcal J}_d^{\theta}|$.  The
following result is a generalition of Theorem \ref{thm2-7} providing
the $h$-principle:
\begin{Theorem}\label{thm2-7-1}
  The map $$j^3_{\pi}(\theta) : |{\mathcal J}_d^{\theta}|\la
  |h{\mathcal J}_d^{\theta}|$$ is weak homotopy equivalence. \hfill
  $\Box$
\end{Theorem}
To describe the homotopy type of the moduli space $\Omega |{\mathcal
  J}^{\theta}_d|$, we consider the bundle diagram:
$$
\begin{diagram}
\setlength{\dgARROWLENGTH}{1.2em}
\node{V^{\theta,\perp}_{d,n}}
     \arrow[2]{e}
     \arrow{s}
\node[2]{V^{\perp}_{d,n}}
     \arrow{s}
\\
\node{G^{\theta,\gmf}(d,n)}
     \arrow[2]{e,t}{\theta_{d,n}}
\node[2]{G^{\gmf}(d,n)}
\end{diagram} 
$$
where $V^{\theta,\perp}_{d,n}$ is a pull-back of
$V^{\perp}_{d,n}$. Similarly, let
$V^{\theta}_{d,n}\to G^{\theta,\gmf}(d,n)$ be the
pull-back of the bundle $V_{d,n}\to G^{\gmf}(d,n)$.  The Thom space of
the bundle $V^{\theta,\perp}_{d,n}$ gives rise to the spectrum
$MT^{\theta,\gmf}(d)$ which in degrees $(d+n)$ is
\begin{equation}\label{eq8b-new}
MT^{\theta,\gmf}(d)_{d+n} = \Th(V^{\theta,\perp}_{d,n}). 
\end{equation}
We have the following version of Theorem \ref{thm3-1}:
\begin{Theorem}\label{thm3-1-1}
There is weak homotopy equivalence
$$
\begin{array}{r}
\ \ \ \ \ \ \ \ \ \ \ \ \ \ \ \ \ \ \ \ \ \ \ \ \ \ \ \ \ \ \ \ \ \ 
\Omega|h{\mathcal J}_d^{\theta}| \cong \Omega^{\infty}MT^{\theta,\gmf}(d) . \ \ \ \ 
\ \ \ \ \ \ \ \ \ \  \  \ \ \ \ 
\ \ \ \ \ \ \ \ \ \  \  \ \ \ \ 
\ \ \ \ \ \ \ \ \ \  
\Box
\end{array}
$$
\end{Theorem}
We next examine the homotopy type of
$\Omega^{\infty}MT^{\theta,\gmf}(d)$ in terms of the corresponding
singular sets $\Sigma^{\theta,\gmf}(d,n)\subset
G^{\theta,\gmf}(d,n)$. Define the $\theta$-Grassmannian
$G^{\theta}(d,n)$ as the pull-back:
\begin{equation}\label{theta-grassm}
\begin{diagram}
\setlength{\dgARROWLENGTH}{1.2em}
\node{G^{\theta}(d,n)}
     \arrow[2]{e}
     \arrow{s,l}{\theta_n}
\node[2]{B}
     \arrow{s,l}{\theta}
\\
\node{G(d,n)}
     \arrow[2]{e,t}{i^o_{n}}
\node[2]{BO(d)}
\end{diagram} 
\end{equation}
where $i^o_n : G(d,n)\rh BO(n)$ is the canonical embedding,
$$
G^{\theta}(d,n)= \{ (V,b) \ | \ i^o_n(V)=\theta(b) \ \}\subset
G(d,n)\times BO(d) .
$$
Then it is easy to identify 
$G^{\theta,\gmf}(d,n)$ with the subspace
$$
G^{\theta,\gmf}(d,n) = \{ \ (V,b,f) \ | \ (V,b)\in G^{\theta}(d,n), \
f\in J^3_{\gmf}(V) \ \}.
$$
Let $\Sigma^{\theta,\gmf}(d,n)$ be the singular set of triples
$(V,b,f)$, where $f: V\to \R$ has vanishing linear term. The
non-sigular subspace $G^{\theta,\gmf}(d,n)\setminus
\Sigma^{\theta,\gmf}(d,n)$ is the set of triples
$(V,b,f)$ with $f=\ell+q+r$, where the linear part $\ell\neq 0$. By
analogy with the space $\hat G^{\gmf}(d,n)$ from Lemma \ref{lem-new1},
we define
$$
\hat G^{\theta,\gmf}(d,n) = \{ \ (V,b,f) \ | \ f=\ell+q+r, \ |\ell|=1, 
\ q = 0, \ r = 0 \ \}
$$
and notice that the space $G^{\theta,\gmf}(d,n)\setminus
\Sigma^{\theta,\gmf}(d,n)$ retracts to $\hat G^{\theta,\gmf}(d,n)$. 
By construction, the space $\hat G^{\theta,\gmf}(d,n)$ is the
pull-back in the diagram:
$$
\begin{diagram}
\setlength{\dgARROWLENGTH}{1.2em}
\node{\hat G^{\theta,\gmf}(d,n)}
     \arrow[2]{e}
     \arrow{s}
\node[2]{\hat G^{\gmf}(d,n)}
     \arrow{s,l}{g_n}
\\
\node{G^{\theta}(d,n)}
     \arrow[2]{e,t}{\theta_n}
\node[2]{G(d,n)}
\end{diagram} 
$$
where $\theta_n: G^{\theta}(d,n)\to G(d,n)$ is from
(\ref{theta-grassm}) and $g_n$ is a composition of the 
forgetting map and the inclusion:
$$
g_n: \hat G^{\gmf}(d,n) \rh G^{\gmf}(d,n) \to  G(d,n)
$$
Similarly to the above case, the normal bundle of the inclusion
$\Sigma^{\theta,\gmf}(d,n)\rh G^{\theta,\gmf}(d,n)$ coincides with
$(V^{\theta}_{d,n})^*\cong V^{\theta}_{d,n}$ restricted to
$\Sigma^{\theta,\gmf}(d,n)$, and  again the inclusion
$$
(D(V^{\theta}_{d,n}),S(V^{\theta}_{d,n}))\longrightarrow
(G^{\theta,\gmf}(d,n), G^{\theta,\gmf}(d,n)\setminus
\Sigma^{\theta,\gmf}(d,n))
$$
is an excision map. Let $j_{\theta}: G^{\theta,\gmf}(d,n)\setminus
\Sigma^{\theta,\gmf}(d,n) \to G^{\theta,\gmf}(d,n)$ be the inclusion.
This leads to the cofibration
\begin{equation}\label{eq1-theta}
\Th(j_{\theta}^*V^{\theta,\perp}_{d,n})\to \Th(V^{\theta,\perp}_{d,n})\to 
\Th((V^{\theta,\perp}_{d,n}\oplus
V^{\theta}_{d,n})_{\Sigma^{\theta,\gmf}(d,n)})
\end{equation}
and to the cofibration of spectra
$$
\Sigma^{-1}MT^{\theta}(d-1)\to MT^{\theta,\gmf}(d)\to
\Sigma^{\infty}(\Sigma^{\theta,\gmf}(d,\infty)_+)
$$
with corresponding homotopy fibration of infinite loop spaces:
$$
\Omega^{\infty}\Sigma^{-1}MT^{\theta}(d-1)\to \Omega^{\infty}MT^{\theta,\gmf}(d)\to
\Omega^{\infty}\Sigma^{\infty}(\Sigma^{\theta,\gmf}(d,\infty)_+).
$$
We denote by $\Cob^{\theta}$ the corresponding cobordism category of
manifolds equipped with tangential structure $\theta$ and by
$\Cob^{\theta,\gmf}$ the category with the condition that each
morphism be equipped with a generalized Morse function as above.
\begin{Corollary}\label{cor1-new-theta}
There is weak homotopy equivalence
$$
B\Cob_d^{\theta,\gmf} \simeq \ \Omega^{\infty-1}MT^{\theta,\gmf}(d),
$$
and the forgetful map $B\Cob_d^{\theta,\gmf}\to B\Cob_d^{\theta}$
induces isomorphism
$$
\Omega^{\theta,\gmf}_d=\pi_0 B\Cob_d^{\theta,\gmf}\cong \pi_0 B\Cob_d=\Omega_d^{\theta},
$$
where $\Omega^{\theta,\gmf}_d$ and $\Omega_d^{\theta}$ are
corresponding cobordism groups.
\end{Corollary}

\begin{small}

\end{small}
\end{document}